\documentclass[11pt,a4paper]{amsart}
\usepackage{amssymb}
\usepackage{amsmath}
\usepackage{graphicx}
\usepackage{color}
\usepackage[utf8x]{inputenc} 	
\usepackage{lmodern}     
\usepackage[T1]{fontenc} 
\usepackage{hyperref}
\newcommand{\be}{\begin{equation}}
\newcommand{\ee}{\end{equation}}
\newcommand{\bea}{\begin{eqnarray}}
\newcommand{\eea}{\end{eqnarray}}

\newcommand{\poly}{\mathbb{P}}
\newcommand{\proj}{\mathcal{P}}
\newcommand{\LL}{\mathcal{L}}
\newcommand{\R}{\mathbb{R}}

\newcommand{\f}{\hat f}
\newcommand{\g}{\hat g}
\newcommand{\Ball}{\mathcal{B}}
\newcommand{\Q}{\mathcal{Q}}
\newcommand{\QL}{\Q^\lambda}
\newcommand{\B}{\hat B}
\newcommand{\bb}{\hat \beta}

\DeclareMathOperator{\Span}{span}

\newtheorem{proposition}{Proposition}
\newtheorem{theorem}{Theorem}
\newtheorem{lemma}{Lemma}
\title{On steady-state preserving spectral methods for homogeneous Boltzmann equations}
\author[F. Filbet]{Francis Filbet}
\address{Francis Filbet \\
Université Lyon 1 \& Inria \\
Institut Camille Jordan \\ 
43 boulevard 11 novembre 1918 \\
69622, Villeurbanne cedex \\
France.}
\email{filbet@math.univ-lyon1.fr}
\author[L. Pareschi]{Lorenzo Pareschi}
\address{Lorenzo Pareschi \\
Mathematics and Computer Science Department \\
University of Ferrara \\
Italy}
\email{lorenzo.pareschi@unife.it}
\author[Th. Rey]{Thomas Rey}
\address{Thomas Rey \\
Center of Scientific Computation and Mathematical Modeling (CSCAMM) \\
The University of Maryland \\
College Park, MD, 20742-4015\\
USA}
\email{trey@cscamm.umd.edu}

\begin{document}

\maketitle

\noindent \emph{\textbf{Abstract.}}
In this note, we present a general way to construct spectral methods for the collision operator of the Boltzmann equation which preserves exactly the Maxwellian steady-state of the system. 
We show that the resulting method is able to approximate with spectral accuracy the solution uniformly in time.  

\bigskip

\noindent \emph{\textbf{Résumé.}} %
Dans cette note, nous présentons une construction générale de méthodes spectrales pour l'opérateur de collision de l'équation de Boltzmann permettant de préserver exactement les états stationnaires Maxwellien de ce type d'équations.
Cette nouvelle approche est basée sur une décomposition de type ``micro-macro'' de la solution de l'équation, tout en restant très proche d'une méthode spectrale plus classique.
Nous montrons que les méthodes obtenues sont capables d'approcher avec une précision spectrale, uniformément en temps, la solution de l'équation considérée, et nous présentons leur efficacité dans un test numérique.

\tableofcontents

\section{Introduction}
Spectral methods for the numerical solution of the homogeneous Boltzmann equation have been proposed originally in~\cite{PP96} and \cite{PaRu:spec:00}. Recently their stability properties were studied in~\cite{FM:11}. 
Related approaches, based on the use of the Fourier transform have been introduced by~\cite{BoRj:HS:99} and~\cite{GT:09}.

In classical kinetic theory the behavior of a
rarefied gas of particles is described by the Boltzmann
equation~\cite{CIP:94}
 \begin{equation*}
 \frac{\partial f}{\partial t} + v \cdot \nabla_x f = \Q(f,f)
 \end{equation*}
where $f(t,x,v)$, $x \in \R^{d}$, $v \in \R^{d}$ ($d \ge 2$), is the
time-dependent particle distribution function in the phase space and
the collision operator $\Q$ is defined by
 \begin{equation}\label{eq:Q}
 \Q (f,f)(v) = \int_{\R^{d}}
 \int_{\mathbb{S}^{d-1}}  B(\cos \theta,|v-v_*|) \,
 \left[ f'_* f' - f_* f \right] \, d\omega \, dv_*.
 \end{equation}
Time and position act only as parameters in $\Q$ and therefore
will be omitted in its description.
In~\eqref{eq:Q} we used the shorthands $f = f(v)$, $f_* = f(v_*)$,
$f ^{'} = f(v')$, $f_* ^{'} = f(v_* ^{'})$. The velocities of the
colliding pairs $(v,v_*)$ and $(v',v'_*)$ are related by
 \begin{equation*}
 v' = \frac{v+v_*}{2} + \frac{|v-v_*|}{2} \sigma, \qquad
 v'_* = \frac{v+v^*}{2} - \frac{|v-v_*|}{2} \sigma\nonumber.
 \end{equation*}
The collision kernel $B$ is a non-negative function which only
depends on $|v-v_*|$ and $\cos \theta = ((v-v_*)/|v-v_*|)\cdot
\omega$.

In this note, we shall assume that the system considered is space homogeneous
and then the particle distribution function $f$ will only depends on the time 
$t$ and the velocity $v$. The Boltzmann equation then becomes
 \begin{equation}
   \label{eq:BE}
   \frac{\partial f}{\partial t}  = \Q(f,f).
 \end{equation}
Boltzmann's collision operator has the fundamental properties of
conserving mass, momentum and energy
 \begin{equation*}
 \int_{v\in{\R}^{d}} \Q(f,f) \, \phi(v)\,dv = 0, \qquad
 \phi(v)=1,v,|v|^2 
 \end{equation*}
and satisfies the well-known Boltzmann's $H$-theorem
 \begin{equation*} 
 - \frac{d}{dt} \int_{v\in{\R}^{d}} f \log f \, dv = - \int_{{\R}^{d}} \Q(f,f)\log(f) \, dv \geq 0.
 \end{equation*}
The functional $- \int f \log f$ is the entropy of the solution.
Boltzmann $H$-theorem implies that any equilibrium distribution
function has the form of a locally Maxwellian distribution
 \begin{equation*}
 M(\rho,u,T)(v)=\frac{\rho}{(2\pi T)^{{d}/2}}
 \exp \left( - \frac{\vert u - v \vert^2} {2T} \right), 
 \end{equation*}
where $\rho,\,u,\,T$ are the density, mean velocity
and temperature of the gas:
 \begin{equation*}
 \rho = \int_{{\R}^d}f(v)dv, \quad u =
 \frac{1}{\rho}\int_{{\R}^d}vf(v)dv, \quad T = \frac{1}{d\rho}
 \int_{{\R}^d}\vert u - v \vert^2f(v)dv.
 \end{equation*}

A major drawback of the classical spectral method is the lack of exact
conservations and, as a consequence the incapacity of the scheme to preserve the Maxwellian steady states of the system. In this paper we overcome this drawback thanks to a new reformulation of the method which permits to preserve the spectral accuracy and to capture the long time behavior of the system.

\section{A steady-state preserving spectral method}

Although the method applies in principle to any collisional kinetic equation which possesses a local Maxwellian-type equilibrium, here we describe the method in the case of the Boltzmann operator.
\subsection{A micro-macro decomposition} 
Let us start with the decomposition 
\be
f=M+g,
\label{eqs5:micmac}
\ee  
with $M$ the local Maxwellian equilibrium and $g$ such that $\int_{\R^d} g\,\phi\,dv=0$, $\phi=1,v,|v|^2$. When
inserted into a Boltzmann-type collision operator, the decomposition \eqref{eqs5:micmac} gives
\be
\Q(f,f)=\LL(M,g)+\Q(g,g),
\label{eq:decom}
\ee
where $\LL(M,g)=\Q(g,M)+\Q(M,g)$ is a linear operator and we used the fact that 
\be
\Q(M,M)=0.
\label{eq:steady}
\ee
There are two major features in the decomposition (\ref{eq:decom}): 
\begin{enumerate}
\item it embeds the identity (\ref{eq:steady}); 
\item the steady state of (\ref{eq:decom}) is given by $g=0$.
\end{enumerate} 
To illustrate the method let us consider now the space homogenous equation \eqref{eq:BE} that we rewrite using the micro-macro decomposition as
\be
\left\{
\begin{aligned}
\frac{\partial g}{\partial t} &= \LL(M,g)+\Q(g,g),\\
f&=M+g.
\end{aligned}
\label{eq:micmac}
\right.
\ee
\subsection{Derivation of the spectral method}
To simplify notations we derived the spectral method in the classical setting introduced in \cite{PaRu:spec:00}, similarly it can be extended to the representation used in \cite{MP:note:04} for the derivation of fast algorithms. Thus, we perform the usual periodization in a bounded domain of the operators $\LL$ and $\Q$ and denote by $\LL^\lambda$ and $\Q^\lambda$ the operators with cut-off on the relative velocity on $\Ball_0(2\lambda\pi)$. 

Let us first set up the mathematical framework of our analysis. For any
$t \geq 0$, $f_N(v,t)$ is a trigonometric polynomial of degree $N$ in
$v$, i.e. $f_N(t) \in \poly^N$ where
\[
\poly^N = \Span\left\{e^{ik\cdot v}\,|\, -N \leq k_j \leq N,\, j=1,\ldots,d
\right\}.
\]
Moreover, let $\proj_N : L^2([-\pi,\pi]^3) \rightarrow \poly^N$ be the
orthogonal projection upon $\poly^N$ in the inner product of
$L^2([-\pi,\pi]^3)$ 
\[
<f-\proj_N f,\phi>=0,\qquad \forall\,\, \phi\,\in\,\poly^N.
\]
We denote the $L^2$-norm by
\[
||f||_2 = (< f, f>)^{1/2}.
\]
With this definition $\proj_N f=f_N$, where $f_N$ is the truncated
Fourier series of $f$ given by 
\begin{equation*}
f_N(v) = \sum_{k=-N}^N \f_k e^{i k \cdot v},
\end{equation*}
with 
\begin{equation*}
\f_k = \frac{1}{(2\pi)^d}\int_{[-\pi,\pi]^d} f(v)
e^{-i k \cdot v }\,dv.
\end{equation*}
We then write the Fourier-Galerkin approximation of the micro-macro equation \eqref{eq:micmac} 
\be
\left\{
\begin{aligned}
\frac{\partial g_N}{\partial t}&=
\LL_N^\lambda(M_N,g_N)+Q_N^\lambda(g_N,g_N),\\
f_N&=M_N+g_N,
\end{aligned}
\label{eq:specc}
\right.
\ee
where $M_N=\proj_N M$, $g_N=\proj_N g$,
$\LL_N^\lambda(M_N,g_N)=\proj_N \LL^\lambda(M_N,g_N)$ and  $\Q_N^\lambda(g_N,g_N)=\proj_N \Q^{\lambda}(g_N,g_N)$. More precisely, using multi-index notations in the above expressions we have 
\begin{equation}
Q_N^\lambda(g_N,g_N) = \sum_{k=-N}^N \left(\sum_{{l,m=-N} \atop {l+m=k}}^N \g_l\,\g_m
\bb(l,m)\right)e^{i k \cdot v},\quad k=-N,\ldots,N,
\label{eq:CF1}
\end{equation}
where the {Boltzmann kernel modes} $\bb(l,m)=\B(l,m)-\B(m,m)$ are given by
\begin{equation}
\B(l,m) = \int_{\Ball_0(2\lambda\pi)}\int_{\mathbb{S}^{d-1}} 
B(\cos \theta,|q|) e^{-i(l\cdot q^++m\cdot q^-)}\,d\omega\,dq. \label{eq:KM}
\end{equation}
In this last identity, $q=v-v_*$ is the relative velocity and $q^+, q^-$ are the following parametrizations of the post-collisional velocities:
  \[
    q^+ = \frac12\left (q+|q|\omega\right ), \quad q^- = \frac12\left (q-|q|\omega\right ).
  \]

It is immediate to show that 
\begin{proposition}
The function $g_N \equiv 0$ is an admissible local equilibrium of the scheme (\ref{eq:specc}) and therefore $f_N=M_N$ is a local equilibrium state. 
\end{proposition}

\section{Spectral accuracy}
In this section we show that (\ref{eq:specc}) is a spectrally accurate approximation to (\ref{eq:micmac}) provided that $M_N$ is a spectrally accurate approximation of $M$. This is clearly guaranteed if initially the support is large enough. Note that due to space homogeneity, $M$ does not change in time and so does $M_N$. It is interesting to observe that the only difference between scheme (\ref{eq:specc}) and the usual spectral method developed on the original formulation 
\be
\begin{aligned}
\frac{\partial f_N}{\partial t}&=\frac{\partial g_N}{\partial t}\\
&=
\Q_N^\lambda(f_N,f_N),\\
&= \LL_N^\lambda(M_N,g_N)+\Q_N^\lambda(g_N,g_N)+\Q_N^\lambda(M_N,M_N)
\end{aligned}
\ee
is due to the constant (in time) term
\be
\Q_N^\lambda(M_N,M_N) \neq 0,
\label{eq:nt}
\ee
which, as we will prove in the sequel, is spectrally small and is not present in (\ref{eq:specc}). 
 
If $f\in H^r_p ([-\pi,\pi]^d)$, where $r \geq 0$ is an integer and
$H^r_p ([-\pi,\pi]^d)$ is the subspace of the Sobolev space
$H^r([-\pi,\pi]^d)$, which consists of periodic functions, we have the following estimate for the spectral accuracy of the classical spectral method \cite{PaRu:spec:00}. 

\begin{theorem}
Let $f \in H_p^r([-\pi,\pi]^d)$, $r\geq 0$ then
\begin{equation}
||\QL(f,f)-\QL_N(f_N,f_N)||_2 \leq \frac{C}{N^r} \left(||f||_{H^r_p} +
||\QL(f_N,f_N)||_{H^r_p}\right),
\end{equation}
\label{th:s5s}
\end{theorem}

Therefore we have the following result for (\ref{eq:nt})
\begin{lemma}
Let $M_N=\proj M$ where $M$ is a local Maxwellian. Then
\begin{equation}
||\QL_N(M_N,M_N)||_2 \leq \frac{C}{N^r} \left(||M||_{H^r_p} +
||\QL(M_N,M_N)||_{H^r_p}\right),
\end{equation}
\label{th:MM}
\end{lemma}
Merging the two results we obtain the spectral accuracy of the new steady state formulation
\begin{theorem}
Let $f \in H_p^r([-\pi,\pi]^3)$, $r\geq 0$ then
\begin{eqnarray}
\nonumber
||\QL(f,f)-\LL_N^\lambda(M_N,g_N)-Q_N^\lambda(g_N,g_N)||_2 &\leq& \frac{C}{N^r} \left(||f||_{H^r_p} + ||M||_{H^r_p}\right.\\[-.25cm]
\\[-.25cm]
\nonumber
&&\left. +
||\QL(f_N,f_N)||_{H^r_p}+||\QL(M_N,M_N)||_{H^r_p}\right).
\end{eqnarray}
\label{th:news}
\end{theorem}  

\section{A numerical example}

We present in this Section a numerical example of our method for the space
homogeneous Boltzmann equation in dimension $2$, with Maxwell molecules:
\[
  B(\cos \theta,|v-v_*|) = \frac{1}{2 \pi}.
\] 
We compare the classical spectral method to the new steady-state preserving one.
For this, we use an an exact solution of the homogeneous Boltzmann equation,
the so called Bobylev-Krook-Wu solution \cite{Bobylev:75,KrookWu:1977}. 
It is given by
\begin{equation*}
f_{BKW}(t,v) = \frac{\exp(-v^2/2S)}{2\pi S^2} \,\left[2\,S-1+\frac{1-S}{2 \,S}\,v^2 \right]
\end{equation*}
with $S = S(t) = 1-\exp(-t/8)/2$. 

For the resolution of the Boltzmann equation \eqref{eq:BE}, we use the fast spectral method of \cite{FMP:06}, with $N = 32$ half-modes in each direction of the box $[-V,V]^2$ for $V = 8$. We take $M = 8$ angular discretizations. Here, we insists that both the classical and steady-states preserving methods rely on the same numerical algorithm. The only difference between the two is the presence of the constant in time term \eqref{eq:nt}.
In particular, both methods have the same computational cost, namely $\mathcal O\left (M N \log_2 N\right )$.

\begin{figure}[!ht]
  \begin{center}
  \includegraphics[scale=1]{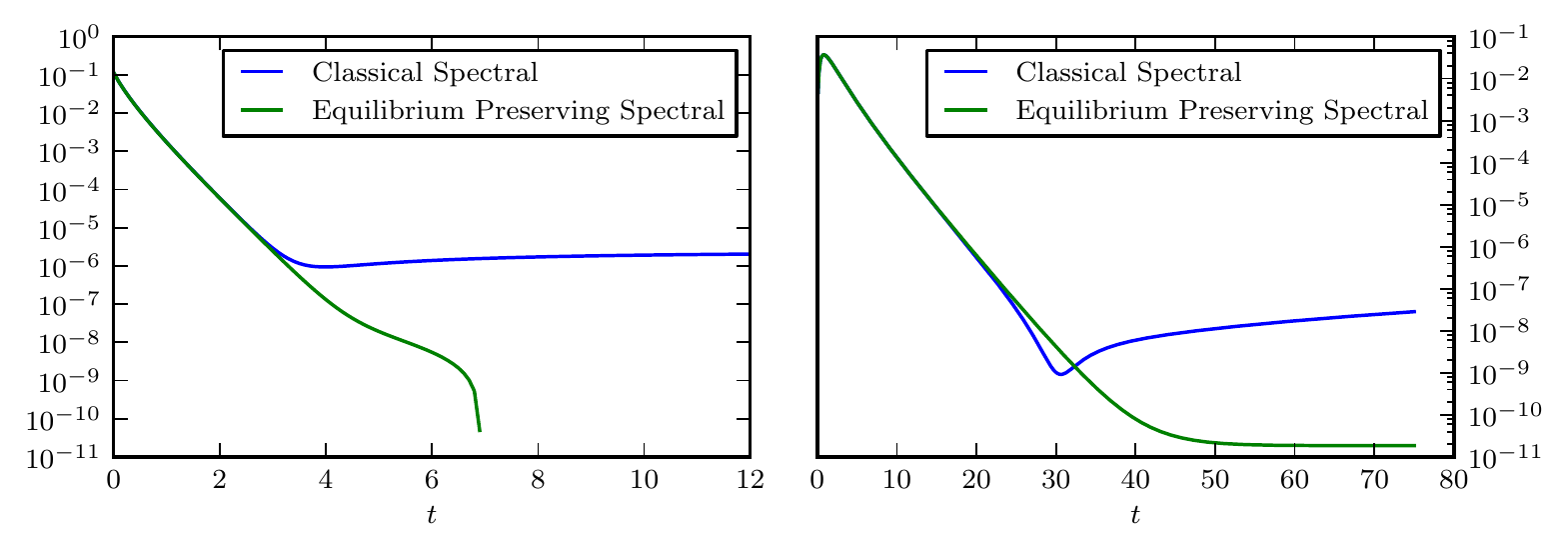}
  \caption{Comparison of the classical spectral method with the steady-state preserving one: Relative entropy $\mathcal H(f|M)(t)$ (left) and $L^2$ error $\|f(t) - f_{BKW}(t)\|_2$ (right).}
  \label{fig:Evolution}
  \end{center}
\end{figure}

Figure \ref{fig:Evolution} presents a comparison between both methods for the relative entropy of the solution $f$ with respect to the global Maxwellian $M$
  \[
    \mathcal H(f|M)(t) := \int_{\R^d} f(t,v) \log \left (\frac{f(t,v)}{M(v)} \right ),
  \]
and the absolute $L^2$ error between the numerical solution $f(t,v)$ and the exact one $f_{BKW}(t,v)$.
We observe that for both quantities, the behavior of the steady-state preserving method is better than the classical one. In particular, the new method achieves a nice monotonous decay of the relative entropy, without the large time increase of the classical spectral method. This is due to the fact that the equilibrium of this latter methods are constants \cite{FM:11}.
The behavior of the $L^2$ error is also monotonous for the steady-state preserving method, which is not the case for the classical one.

\begin{figure}[!ht]
  \begin{center}
  \includegraphics[scale=1]{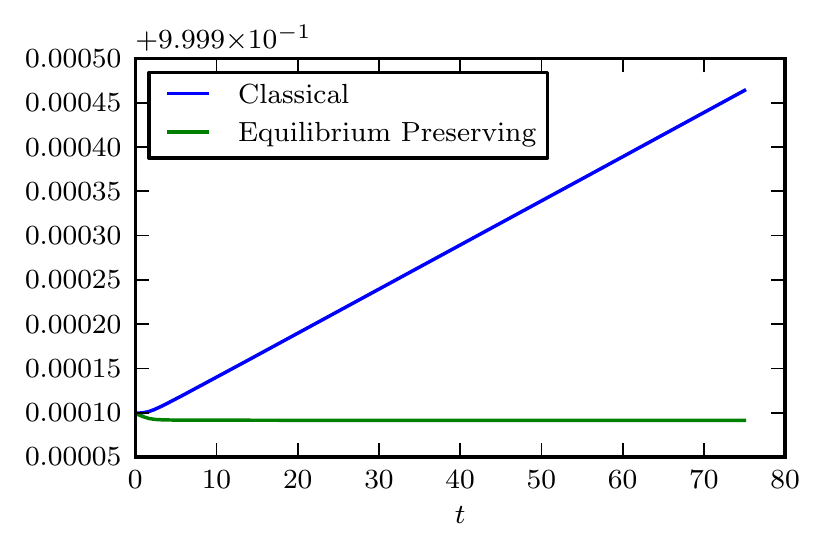}
  \caption{Comparison of the classical spectral method with the steady-state preserving one: Evolution of the temperature.}
  \label{fig:Temperature}
  \end{center}
\end{figure}

We then present in Figure 2 a comparison between the temperatures of the solutions obtained with the classical and steady-state preserving methods. As expected with such a few Fourier modes, the temperature of the classical spectral method is not preserved, but grows almost linearly because of the convergence toward the discrete equilibrium, a constant on the cube.
In contrast, it seems that the steady-state preserving spectral methods improve greatly this behavior: although not perfectly conserved, the temperature of the solution decays by roughly $0.1\%$ and then stabilizes on a constant. 

Finally, as a last observation, we noticed that even if the positivity of the solution is not preserved by the new method, the number of nonpositive cells seems to be lower at a fixed time for the steady-state preserving method, compared to the classical one.
We recall that it was already proved in \cite{FM:11} that the number of nonpositive cells is asymptotically $0$ for the classical spectral method:

\begin{theorem}[Filbet, Mouhot, \cite{FM:11}]
  Let $f_N$ be the unique, global solution to the truncated, Fourier transformed Boltzmann equation, associated to a nonnegative initial datum $f_0 \in H^k_{p}([-\pi,\pi]^d)$ with $k > d/2$. Then the mass of the negative values of $f_N$ can be made uniformly in time small as the number of Fourier modes $N \to \infty$. 
\end{theorem}

	\section*{Acknowledgment}
	  The research of FF  is partially supported by the Inria EPI Kaliffe and European Research Council ERC Starting Grant 2009, project 239983-NuSiKiMo.
	  The research of TR is granted by the NSF Grants  \#1008397 and \#1107444 (KI-Net) and ONR grant \#000141210318.
	  Part of this work was conducted during a KI-Net (NSF \#1107444) meeting in North Carolina State University.

\end{document}